\documentclass[11pt,reqno]{amsart}
\usepackage{amsfonts}
\usepackage{fancyhdr}
\usepackage{times}
\usepackage[T1]{fontenc}
\usepackage{mathrsfs}
\usepackage{latexsym}
\usepackage[dvips]{graphics}
\usepackage{epsfig}
\usepackage{hyperref, amsmath, amsthm, amsfonts, amscd, flafter,epsf}
\usepackage{amsmath,amsfonts,amsthm,amssymb,amscd}
\input amssym.def
\input amssym.tex
\usepackage{color}
\usepackage[parfill]{parskip}
 \usepackage{mathrsfs} 
\usepackage{geometry}                
\usepackage{epstopdf}

\usepackage{verbatim}

\newtheorem{thm}{Theorem}

\newtheorem{lem}[thm]{Lemma}
\newtheorem{thm*}{Theorem}
\newtheorem*{con*}{Conjecture}
\newtheorem*{lem*}{Lemma}

\def\g{\gamma}
\def\b{\beta'-1/2}

\begin{document}
 \baselineskip=17pt
\hbox{}
\medskip
\title{Conditional estimates on small distances between \\ordinates of zeros of $\zeta(s)$ and $\zeta'(s)$}
\author{Fan Ge}

\email{fan.ge@rochester.edu}

\address{Department of Mathematics, University of Rochester, Rochester, NY}

\thanks{This work was partially supported by NSF grant DMS-1200582.}

\maketitle

\begin{abstract}
Let $\beta'+i\g'$ be a zero of $\zeta'(s)$.
In \cite{GYi} Garaev and  Y{\i}ld{\i}r{\i}m proved that
there is a zero $\beta+i\g$ of $\zeta(s)$ with
$
\g'-\g \ll \sqrt{|\beta'-1/2|}
$.
Assuming RH, we improve this bound by saving a factor $\sqrt{\log\log\g'}$.
\end{abstract}

\section{introduction}

The distribution of zeros of the Riemann zeta-function $\zeta(s)$ is
closely connected to that of zeros of $\zeta'(s)$. As just one
illustration we cite A. Speiser's~\cite{Spe} theorem that the
Riemann Hypothesis (RH) is equivalent to the nonexistence of
non-real zeros of $\zeta'(s)$ in the half-plane $\Re s<1/2$.

Let $\rho'=\beta'+i\g'$ be a zero of $\zeta'(s)$, and let
$\rho_c=\rho_c(\rho')=\beta_c+i\g_c$ be a zero of $\zeta(s)$ with
smallest $|\g'-\g_c|$ (if there is more than one such zero, take any
of them). M. Z. Garaev and C. Y. Y{\i}ld{\i}r{\i}m ~\cite{GYi}
showed that
$$
\g'-\g_c \ll \sqrt{|\beta'-1/2|}.
$$
Their result is unconditional. Our purpose here is to obtain a
conditional improvement.
\begin{thm}\label{thm}
Assume RH. We have
\begin{align}\label{eq thm}
\g'-\g_c \ll \sqrt{\frac{\beta'-1/2}{\log\log\g'}}
\end{align}
for $\beta'-1/2\le1/\log\log\g'$. Here the implied constant is absolute,
and for $\g'$ sufficiently large we may take the implied constant to be 2.16.
\end{thm}

\emph{Remark 1.}  Note that on RH we trivially have
\begin{align}
\label{eq trivial bound}
\g' - \g_c \ll \frac{1}{\log\log \g'}\ .
\end{align}
Combining this with our Theorem \ref{thm}, we see that on RH
\begin{align*}
\g'-\g_c \ll \min \bigg\{\sqrt{\frac{\beta'-1/2}{\log\log\g'}},\
\frac{1}{\log\log \g'}\bigg\}.
\end{align*}
 The inequality \eqref{eq trivial bound} follows
from the well-known fact that on RH, the largest gap between
consecutive zeros of $\zeta(s)$ up to height $T$ is $\ll 1/\log\log
T$ (see \cite{Tit}, for example).

\emph{Remark 2.} In \cite{FGH} D. W. Farmer, S. M. Gonek and C. P.
Hughes conjectured that
$$
\limsup_{t\rightarrow \infty} \frac{S(t)}{\sqrt{\log t \log\log t}} = \frac{1}{\pi \sqrt{2}}.
$$
Assuming this as well as RH, one can show (by the same proof as that
of Theorem \ref{thm}) that
$$ \g'
- \g_c \ll \sqrt{\beta'-1/2}\ \bigg(\frac{\log\log
\g'}{\log \g'}\bigg)^{1/4}
$$
for $\beta'-1/2\ll \sqrt{\log\log \g'/\log \g'}$.

\emph{Remark 3.} There are multiple ways to prove results like
Theorem 1. For example, one can start with Lemma~\ref{lem sum of h}
below, split the sum into three parts (according to
$|\g-\g'|\le1/\log\log \g'$, $1/\log\log \g'<|\g-\g'|\le1$ or
$|\g-\g'|\ge1$), and estimate each part separately. This will give a
slightly weaker result than Theorem 1. The proof we present in this
paper follows another clue, which we think is more inspiring and
more likely to be modified. For example, with a little more care it
is possible to show that (on RH) for $\beta'-1/2\le1/\log\log\g'$
there are $ \gg (\b)\log\g' $ zero(s) of $\zeta(s)$ lie in $
\Big[\,\g'-C\sqrt{\frac{\beta'-1/2}{\log\log\g'}},\
\g'+C\sqrt{\frac{\beta'-1/2}{\log\log\g'}}\ \Big] $ for some
constant $C$.

\section{lemmas}

\begin{lem}\label{lem sum of h}
Assume RH. If $\beta'>1/2$, then we have
$$
\frac{ \log\g'}{2} = \sum_{\g} \frac{\b}{(\b)^2+(\g'-\g)^2}+O(1).
$$
\end{lem}
See equation (4) in \cite{Sou}.

Let $ N(T)=\sum_{0<\g\leq T} 1 $ be the zero counting function of
$\zeta(s)$. It is well-known (see \cite{Tit}) that
$$N(T)=L(T)+S(T)+E(T),$$ where \begin{align*}L(T)=\frac{1}{2\pi}T\log T - \frac{1+\log 2\pi}{2\pi}T
+\frac{7}{8}, \ \ \ \ \ \
S(T)=\pi^{-1}\arg\zeta(1/2+iT),\end{align*} and $E(T)$ is an error
term. We require the following result.

\begin{lem}\label{lem L+E}
We have $$d(L(u)+E(u))=\bigg(  \frac{1}{2\pi}\log u +O(1)  \bigg)du.   $$
\end{lem}
\proof By the proof of Theorem 9.3 in \cite{Tit} we know that
$$
L(T)+E(T)=1-\frac{T\log\pi}{2\pi} + \frac{1}{\pi}\Im\log\Gamma(1/4+iT/2).
$$
Therefore, we have $$d(L(u)+E(u))=\bigg( \frac{1}{\pi}\frac{d\
\Im\log\Gamma(1/4+iu/2)}{du}  +O(1) \bigg)du.$$ It is
straightforward to compute that
$$ \frac{d\ \Im\log\Gamma(1/4+iu/2)}{du} =
\frac{1}{2}\ \Re\frac{\Gamma'}{\Gamma}(1/4+iu/2).$$ By Stirling's
formula, this is $(\log u) /2 +O(1)$. Hence the result. \qed

\begin{lem}\label{lem E difference}
Let $T>2$ and $T<t_1<t_2<2T$. Then we have
$$E(t_2)-E(t_1)\ll t_2-t_1.$$
\end{lem}
\proof Write $s=1/2+it$. By the proof of Theorem 9.3 in \cite{Tit} we know that
\begin{align*}
\pi E(t)& = \pi (N(t)-S(t)-L(t))\\
    & = \Delta\arg s(s-1) + \Delta\arg \pi^{-s/2} +
    \Delta\arg\Gamma(s/2) +\Delta\arg \zeta(s) \\ & \qquad - \arg\zeta(s) -
    \frac{T}{2}\log t +\frac{1+\log 2\pi}{2}t - \frac{7}{8}\pi\\
    & = \Delta\arg\Gamma(s/2) -\frac{T}{2}\log t
    +\frac{1+\log 2}{2}t + \frac{\pi}{8}.
\end{align*}
It follows that
$$
\pi (E(t_2)-E(t_1)) =  \Delta\arg\Gamma(1/4+it_2/2)-\Delta\arg\Gamma(1/4+it_1/2) - \frac{1}{2}(t_2-t_1)\log T +O(t_2-t_1).
$$
By the mean value theorem of calculus,
 $$\Delta\arg\Gamma(1/4+it_2/2)-\Delta\arg\Gamma(1/4+it_1/2)=(t_2-t_1)\cdot \frac{1}{2}\ \Re\frac{\Gamma'}{\Gamma}(1/4+it_3/2)
 $$
 for some $t_3\in [t_1, t_2]$. But this is
 $$
  \frac{1}{2}(t_2-t_1)\log T
  +O\bigg(\frac{t_2-t_1}{T}\bigg)
  $$ by Stirling's formula.
Hence the result. \qed

\section{proof of the theorem}

It is well-known that $$\zeta'(1/2+i\g')=0 \Longrightarrow
\zeta(1/2+i\g')=0.$$ Therefore, $\beta'=1/2$ implies that
$\g_c=\g'$, in which case \eqref{eq thm} is trivially true. Below we
assume that $1/2<\beta'\le 1/2+1/\log\log \g'$. We may also assume
$\g'>2015$ for convenience.

Define
 $$
 h(t)= h_{\rho'}(t) = \frac{\b}{(\b)^2+(t-\g')^2}\ .
 $$
By Lemma \ref{lem sum of h} we have
$$
\sum_{\gamma} h(\gamma)= \frac{1}{2}\log\g'+O(1).
$$
It is well-known that $\zeta(s)$ has no zero in the region
$$
\sigma>0,\ \ \ \ -14\le t\le 14.
$$
For $t\le-14$, there are $\ll \log |t|$ zeros $1/2+i\g$ of $\zeta(s)$ for which $t-1\le \g \le t$. Thus, it is easy to see that
$$
 \sum_{-\infty<\g\le 14} h(\g) = \sum_{n=14}^{\infty} \ \ \sum_{-n-1<\g\le -n} h(\g) \ll (\b)\cdot \sum_{n=1}^{\infty}\frac{\log
 n}{n^2}\ll 1.$$
It follows that
\begin{align} \label{eq whole integral}
\frac{1}{2}\log\g'+O(1)  = \sum_{\gamma>14} h(\gamma)   =
\int_{14}^{\infty}h(u)d(N(u))                    =
\int_{14}^{\infty}h(u)d(L(u)+E(u)+S(u)).
\end{align}

Next we show that
$$\int_{14}^{\infty}h(u)d(L(u)+E(u)) = \frac{\log\g'}{2}+O(1).
 $$

By Lemma \ref{lem L+E} we have
$$
\int_{14}^{\infty}h(u)d(L(u)+E(u)) = \int_{14}^{\infty}h(u)\bigg(\frac{\log u}{2\pi}+O(1)\bigg)du.
$$
It is clear that $$\bigg(\int_{14}^{\g'/2} + \int_{3\g'/2}^{\infty}\bigg)\bigg( h(u)\Big(\frac{\log u}{2\pi}+O(1)\Big)\bigg)du \ll 1,$$
and that $$\int_{\g'/2}^{3\g'/2}h(u)\bigg(\frac{\log u}{2\pi}+O(1)\bigg)du = \frac{\log\g'}{2\pi}\int_{\g'/2}^{3\g'/2}h(u)du+O(1).$$
Hence, we see that
$$
\int_{14}^{\infty}h(u)d(L(u)+E(u)) =\frac{\log\g'}{2\pi}\int_{\g'/2}^{3\g'/2}h(u)du+O(1).
$$
Now we plainly have
\begin{align*}
\int_{\g'/2}^{3\g'/2}h(u)du  = 2\arctan\bigg(\frac{\g'}{2(\b)}\bigg)
  = \pi + O\bigg(\frac{\b}{\g'}\bigg).
 \end{align*}
 Therefore, we obtain
 $$\int_{14}^{\infty}h(u)d(L(u)+E(u)) = \frac{\log\g'}{2}+O(1).
 $$
This together with (\ref{eq whole integral}) give us
$$
\int_{14}^{\infty}h(u)dS(u)=O(1).
$$

By integration by parts, we see that
$$
 \int_{14}^{\infty}h(u)dS(u) = -\int_{14}^{\infty}h'(u)S(u)du +O(1).
 $$
 It follows that
 \begin{align}\label{eq whole h'S}
 \int_{14}^{\infty}h'(u)S(u)du=O(1).
 \end{align}

Let $p=p(\g')$ be a parameter to be determined later. Split the above integral into three parts:
\begin{align*}
 \int_{14}^{\infty}h'(u)S(u)du  = \bigg[\Big( \int_{14}^{\g'-\frac{\sqrt{\b}}{p}} + \int_{\g'+\frac{\sqrt{\b}}{p}}^{2\g'}\Big)+
 \int_{\g'-\frac{\sqrt{\b}}{p}}^{\g'+\frac{\sqrt{\b}}{p}} + \int_{2\g'}^\infty\bigg]h'(u)S(u)du\ .
 \end{align*}
 We estimate them separately.
First, since
$$
h'(u)=-\frac{2(u-\g')(\b)}{((\b)^2+(u-\g')^2)^2}\ ,
$$
we trivially have
\begin{align}
\label{eq int fourth part}
\int_{2\g'}^\infty h'(u)S(u)du \ll 1/\g'.
\end{align}

Next we consider
$$
\Big( \int_{14}^{\g'-\frac{\sqrt{\b}}{p}} + \int_{\g'+\frac{\sqrt{\b}}{p}}^{2\g'}\Big)h'(u)S(u)du
\ .
$$
It is straightforward to compute that
$$
\int_{-\infty}^\infty   |h'(u)|du = 2\int_{-\infty}^{\g'} h'(u) du =
2 h(u)\Big|_{-\infty}^{\g'}= \frac{2}{\b}\ ,
$$
and that
$$
\int_{\g'-\frac{\sqrt{\b}}{p}}^{\g'+\frac{\sqrt{\b}}{p}}|h'(u)|du = \frac{2}{\b}\cdot\frac{1}{1+(\b)p^2}\ .
$$
Hence, using the bound (see \cite{Tit})
\begin{align*}
|S(T)|\le  \frac{ A\log T}{\log\log T}
\end{align*}
for some absolute positive constant $A$,
we see that
\begin{align}  \label{eq int two parts}
\bigg(\int_{14}^{\g'-\frac{\sqrt{\b}}{p}} + &
\int_{\g'+\frac{\sqrt{\b}}{p}}^{2\g'}\bigg) \big|h'(u)S(u)\big|du \nonumber\\ & \le
\frac{ 2A \log\g'}{\log\log\g'}\bigg(\int_{14}^{\g'-\frac{\sqrt{\b}}{p}} +
\int_{\g'+\frac{\sqrt{\b}}{p}}^{2\g'}\bigg) \big|h'(u)\big|du \nonumber\\ &
\le \frac{2A\log\g'}{\log\log\g'}\bigg(\int_{-\infty}^\infty -
\int_{\g'-\frac{\sqrt{\b}}{p}}^{\g'+\frac{\sqrt{\b}}{p}}\bigg)\big|h'(u)\big|du  \nonumber
\\ &
=  \frac{2A\log\g'}{\log\log\g'} \bigg( \frac{2}{\b} -   \frac{2}{\b}\cdot\frac{1}{1+(\b)p^2}    \nonumber
\bigg)\\
& = \frac{2A\log\g'}{\log\log\g'}\cdot \frac{2p^2}{1+(\b)p^2}\ .
\end{align}

Now we turn to
$$
\int_{\g'-\frac{\sqrt{\b}}{p}}^{\g'+\frac{\sqrt{\b}}{p}} h'(u)S(u)du\ .
$$
Suppose that there is no zero of $\zeta(s)$ on the vertical segment
\begin{align}\label{assump}
\bigg[1/2+i\Big(\g'-\frac{\sqrt{\b}}{p}\Big),\
1/2+i\Big(\g'+\frac{\sqrt{\b}}{p}\Big)\bigg]\ .
\end{align}
Then we
have $ N(t_2)-N(t_1)=0 $ for $ t_1, t_2\in
\Big[\,\g'-\frac{\sqrt{\b}}{p},\g'+\frac{\sqrt{\b}}{p}\ \Big] $. It
follows that
$$
S(t_1)-S(t_2)=L(t_2)-L(t_1)+E(t_2)-E(t_1)=\frac{t_2-t_1}{2\pi}\log\g'+O(t_2-t_1)+E(t_2)-E(t_1).
$$
By Lemma \ref{lem E difference}, this is
\begin{align}\label{eq dif s}
\frac{t_2-t_1}{2\pi}\log\g'+O(t_2-t_1).
\end{align}
Therefore, since
$$
h'(u)=-\frac{2(u-\g')(\b)}{((\b)^2+(u-\g')^2)^2}\ ,
$$
by changing variables we see that
 $$
 \int_{\g'-\frac{\sqrt{\b}}{p}}^{\g'+\frac{\sqrt{\b}}{p}}
 h'(u)S(u)du  = \int_0^{\frac{\sqrt{\b}}{p}}
 \frac{2(\b)v}{((\b)^2+v^2)^2}\bigg(S(\g'-v)-S(\g'+v)\bigg)dv\ .
 $$
 By \eqref{eq dif s}, this is
 $$
 \int_{\g'-\frac{\sqrt{\b}}{p}}^{\g'+\frac{\sqrt{\b}}{p}}
 h'(u)S(u)du  =
 \int_0^{\frac{\sqrt{\b}}{p}}
 \frac{4(\b)v^2}{((\b)^2+v^2)^2}\bigg(\frac{\log
 \g'}{2\pi}+O(1)\bigg)dv\ ,
 $$
 and a straightforward computation turns it into
 $$
 \int_{\g'-\frac{\sqrt{\b}}{p}}^{\g'+\frac{\sqrt{\b}}{p}}
 h'(u)S(u)du  =
  \bigg(\frac{\log
 \g'}{\pi}+O(1)\bigg)\cdot
 \bigg(\frac{-p\sqrt{\b}}{1+(\b)p^2}+\arctan\Big(\frac{1}{p\sqrt{\b}}
\Big)\bigg).
$$
Combining this with (\ref{eq whole h'S}), (\ref{eq int fourth part}) and (\ref{eq int two parts}) we obtain
\begin{align}
\label{eq contradiction}
\bigg(\frac{\log
 \g'}{\pi}+O(1)\bigg)\cdot&
 \bigg(\frac{-p\sqrt{\b}}{1+(\b)p^2}+\arctan\Big(\frac{1}{p\sqrt{\b}}\nonumber
\Big)\bigg)\\ \nonumber & = \int_{\g'-\frac{\sqrt{\b}}{p}}^{\g'+\frac{\sqrt{\b}}{p}}
 h'(u)S(u)du\\ & = \bigg[ \int_{14}^{\infty} - \Big( \int_{14}^{\g'-\frac{\sqrt{\b}}{p}} +\nonumber
 \int_{\g'+\frac{\sqrt{\b}}{p}}^{2\g'}\Big)- \int_{2\g'}^\infty\bigg]h'(u)S(u)du\\ & \le
O(1) +\  \frac{2A\log\g'}{\log\log\g'}\cdot
\frac{2p^2}{1+(\b)p^2} +\ O(1/\g')\ .
\end{align}

We wish to choose $p=c\sqrt{\log\log\g'}$ for some sufficiently
small positive constant $c$ such that
\begin{align}\label{eq contrdct 1}
\frac{-p\sqrt{\b}}{1+(\b)p^2}+\arctan\Big(\frac{1}{p\sqrt{\b}}
\Big) \ge \frac{\pi}{3}\ ,
\end{align}
and that
\begin{align}\label{eq contrdct 2}
\frac{2A}{\log\log\g'}\cdot
\frac{2p^2}{1+(\b)p^2}\le \frac{1}{6}.
\end{align}

We show such $c$ exists. In fact, we clearly have
$$
\frac{2A}{\log\log\g'}\cdot \frac{2p^2}{1+(\b)p^2} = \frac{4Ac^2}{1+c^2(\b)\log\log\g'}\le 4Ac^2.
$$
Next, since $\b\le 1/\log\log\g'$ we have $0<p\sqrt{\b}\le c$. It follows that
$$
\frac{-p\sqrt{\b}}{1+(\b)p^2}+\arctan\Big(\frac{1}{p\sqrt{\b}}
\Big)\ge -c + \arctan(c^{-1}).
$$
Thus, there does exist a small constant $c>0$ such that both (\ref{eq contrdct 1}) and
(\ref{eq contrdct 2}) hold.

Now combining (\ref{eq contradiction}) with (\ref{eq contrdct 1}) and
(\ref{eq contrdct 2}), we obtain
$$
\bigg(\frac{\log
 \g'}{\pi}+O(1)\bigg)\cdot \frac{\pi}{3} \le \frac{\log \g'}{6} + O(1),
 $$
which is clearly a contradiction for large $\g'$.

Hence, the assumption \eqref{assump} must be false. This means there
exists a zero of $\zeta(s)$ on the vertical segment
$$
\bigg[1/2+i\Big(\g'-\frac{\sqrt{\b}}{p}\Big),\
1/2+i\Big(\g'+\frac{\sqrt{\b}}{p}\Big)\bigg]\ . $$ This ends our
proof. \qed

\emph{Note added in proof}:
From the above discussion we see that for any $\epsilon>0$ and $\g'$ sufficiently large (depending on $\epsilon$),
it suffices to choose $c$ such that
\begin{align*}
-c + \arctan(c^{-1}) \ge
\pi\cdot\ 4Ac^2\ +\epsilon.
\end{align*}
By the work of E. Carneiro, V. Chandee and M. B. Milinovich \cite{CCM},
we can take $A=\frac{1}{4}+o(1)$. Therefore, we may choose any positive $c<c_0$ where $c_0$ is the positive root of $\arctan(x^{-1})-x = \pi x^2$,
whose numerical value is $c_0=0.463...$. Thus, for $\g'$ sufficiently large, we may then take the implied constant in \eqref{eq thm} to be $1/0.463\approx 2.16$.

\section*{Acknowledgement}
    The author is indebted to Professor Steve Gonek for very helpful
    conversations. He also thanks the referee for valuable suggestions.


\begin{thebibliography}{99}

\bibitem{CCM} E. Carneiro, V. Chandee and M. B. Milinovich, \textit{Bounding $S(t)$ and $S_1(t)$ on the Riemann hypothesis}, Math. Ann. 356 (2013), 939-968.

\bibitem{FGH} D. W. Farmer, S. M. Gonek and C. P. Hughes, \textit{The maximum size of L-functions}, J. Reine Angew. Math. 609 (2007), 215-236.

\bibitem{GYi} M. Z. Garaev and C. Y. Y{\i}ld{\i}r{\i}m, \textit{On small distances between ordinates of zeros of $\zeta(s)$ and $\zeta'(s)$}, Int. Math. Res. Not. IMRN 2007 (21) (2007), Art. ID rnm091, 14pp.

\bibitem{GG} D. A. Goldston and S. M. Gonek, \textit{A note on $S(t)$ and the zeros of the Riemann zeta-function}, Bull.
Lond. Math. Soc. 39 (2007), 482-486.

\bibitem{LMo} N. Levinson and H. L. Montgomery, \textit{Zeros of the derivatives of the Riemann zeta-function}, Acta Math. 133 (1974), 49-65.

\bibitem{Spe} A. Speiser, \textit{Geometrisches zur Riemannschen Zetafunktion}, Math. Ann. 110 (1934), 514-521.


\bibitem{Sou} K. Soundararajan, \textit{Moments of the Riemann zeta-function}, Ann. of Math. (2), 170(2) 981-993, 2009

\bibitem{Tit} E. C. Titchmarsh, \textit{The theory of the Riemann zeta-function}, 2nd ed., (ed. D. R. Heath-Brown; Oxford Science Publications, Oxford, 1986).



\end{thebibliography}
\end{document}